\documentclass[12pt, a4paper]{article}
\usepackage{amsfonts}
\usepackage{amscd}
\usepackage{graphicx}
\usepackage{amsmath}
\usepackage{caption}
\usepackage{amssymb}
\usepackage{amsthm}

\newtheorem{theorem}{Theorem}

\theoremstyle{definition}

\newtheorem{definition}{Definition}

\begin{document}
\title{\bf On the ellipticity of operators associated with Morse-Smale diffeomorphisms}
\date{}
\author{N.R. Izvarina, A.Yu. Savin}

\maketitle
\abstract{
We consider the operator algebra generated by pseudodifferential operators on a closed smooth surface and shift operator induced by a Morse--Smale diffeomorphism of this surface. Elements in this algebra are considered as operators in the scale of Sobolev spaces and the aim of this paper is to describe how Fredholm property of a given operator depends on the Sobolev smoothness exponent $s$.}

\section*{Introduction}

Given a group action on a manifold, consider the following operator algebra: its elements are linear combinations of shift operators along the orbits of the group action with pseudodifferential operators (below $\psi$DO) on the manifold as coefficients. Such operators were studied by many authors, e.g., see the monographs
\cite{AnLe1,AnLe2,Con1,NaSaSt17} and the references cited there.  

It is important to note that ellipticity (and, hence, the Fredholm property) for operators in this algebra  in general depends on the smoothness exponent $s$ of the Sobolev space   $H^s$, on which our operator is considered. In particular, an operator might be elliptic for some values of $s$ and not elliptic for  other values of $s$.  Hence, there is a natural problem of determining the values of $s$, for which a given operator is elliptic. 

It is known that ellipticity does not depend on $s$ for isometric actions, hence, nonisometric actions should be studied. First steps were undertaken in  \cite{SaSt25,Sav12} 
in the study of operators associated with nonisometric actions. In particular, for the group generated by a dilation
of the sphere, it was shown that a given operator is elliptic for all $s$ in a certain interval (possibly infinite, semi-infinite, or even empty) and is not elliptic for all $s$ outside this interval. Then, in   \cite{Izv1} 
operators associated with a  parabolic diffeomorphism of spheres were studied. Even though this diffeomorphism is nonisometric, it was shown that nonetheless, ellipticity does not depend on $s$ in this case.  

In this paper, we study operators, associated with diffeomorphisms of Morse-Smale type on surfaces.  For operators with constant coefficients, we give explicit conditions for their invertibility in the scale of Sobolev spaces
in terms of their coefficients and invariants of the diffeomorphism at the fixed points. For general operators,
we show that  there exists an interval $I$ (this interval depends on the operator) such that the operator is not elliptic for all $s\notin I$ and it is either elliptic or not elliptic for all $s\in I$ simultaneously.  

This work was supported by RFBR grant Nr.16-01-00373 and also by G-RISC project Nr.~M-2017-b.
We are also grateful to Prof. E.~Schrohe (Leibniz University of Hannover) for his interest to this work.

\section{Statement of the problem}

Given a diffeomorphism $g:M\to M$ of a closed smooth manifold   $M$, we consider operators
equal to finite sums of the form  
\begin{equation}
\label{1}
D=\sum_k{D_kT^k}: \ \  H^s(M) \longrightarrow H^{s-d}(M),\ \ \ s\in\mathbb R,
\end{equation}
where $D_k$ are $\psi$DOs of order $d$ on $M$, $k\in\mathbb Z$, while $T$ is a shift operator associated with  $g: Tu(x)=u(g(x))$. Let us recall the ellipticity conditions for such operators, which guarantee that the operator \eqref{1}
has Fredholm property. As usual, these conditions are formulated in terms of the symbol of operator. The symbol is defined as follows, see \cite{AnLe1,Sav12}. Let us recall that $g$ induces a diffeomorphism
$$
\partial g: T^*M\longrightarrow T^*M
$$ 
of the cotangent bundle of $M$. This diffeomorphism is called the {\em codifferential} of $g$ and is defined as
$\partial g=((dg)^t)^{-1}$, where $dg:TM\to TM$ is the differential of $g$.
The {\em symbol} of $D$ at a point $(x,\xi)\in T^*_0M$ is the operator
\begin{equation}
\label{2}
\begin{array}{ccc}
\sigma(D)(x,\xi):l^2(\mathbb Z, \mu_{x,\xi,s}) & \longrightarrow & l^2(\mathbb Z, \mu_{x,\xi,s-d}) \vspace{4mm}\\
\sigma(D)(x,\xi)w(n) &= &\sum_k \sigma(D_k)(\partial g^n(x,\xi))\mathcal T^k w(n),
\end{array}
\end{equation}
where $\sigma(D_k)$ is the symbol of $D_k$, while $\mathcal Tu(n)=u(n+1)$ is the shift operator.  
Note that the symbol is a family parametrized by $(x,\xi)$ of finite-difference operators with variable coefficients.
These operators act on the spaces $l^2(\mathbb Z, \mu_{x,\xi,s})$ of sequences $\{w(n)\}_{n\in\mathbb{Z}}$ square-summable 
$$
\sum_{n} |w(n)|^2 \mu_{x,\xi,s}(n)<\infty
$$
with respect to the weight $\mu_{x,\xi,s}$ defined as
\begin{equation}
\label{3}
\mu_{x,\xi,s}(n) = \frac{\partial g^{n*}[\mu\sigma(\Delta^s)](x,\xi)}{\mu(x)},
\end{equation}
where $\mu$ is a smooth nonsingular measure on $M$, while $\Delta$ stands for the Laplacian on $M$.
 
Operator \eqref{1} is {\em elliptic} if its symbol is invertible for all $(x,\xi)\in T^*_0M$. It is known, e.g., see \cite{AnLe2,Sav12}, 
that ellipticity implies Fredholm property. Since the symbol depends on the Sobolev smoothness exponent $s$,  ellipticity condition in general depends on $s$. Unfortunately, in the general case this dependence is difficult to describe, because the asymptotic  of the weights as $n\to\infty$ might be quite complicated (depending on the dynamics of $g$). 

\begin{figure}[h!]
	\center{\includegraphics[width=0.5\linewidth]{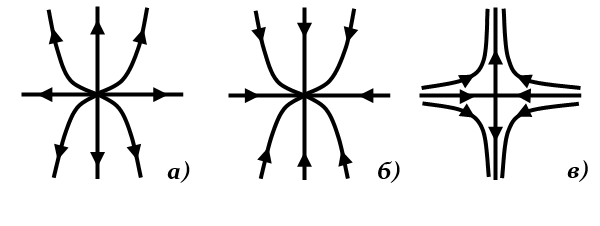}}
	\caption{ Types of fixed points: source \sl{a)}; sink \sl{b)}; saddle \sl{c)}.}
	\label{fig:im1}
\end{figure}

In this paper, we study the case of operators   \eqref{1}, associated with diffeomorphisms of Morse-Smale type on surfaces   (in this situation, the weights in \eqref{3} can be explicitly calculated).
\begin{definition}\label{def1}
We say that a diffeomorphism $g:M\to M$ is of {\em Morse--Smale type}, if it   satisfies the following conditions:
\begin{enumerate}
\item The set $M^g$ of fixed points is finite and there are no periodic points with period $\ge 2$.
\item Given a fixed point $x_0\in M^g$, the eigenvalues of the differential 
$$
dg|_{x=x_0} : \ \ T_{x_0}M \longrightarrow T_{x_0}M
$$
are real, positive, distinct  and not equal to 1.
\end{enumerate} 
\end{definition}
From now on we suppose that $g$ satisfies the conditions in this definition.
Thus, on   $M$, we have a finite number of fixed points of the following types:  source, sink, saddle (see Fig.~\ref{fig:im1}). 
In addition, given a point $x\in M$, the points $g^nx$ have limits as $n\to\pm\infty$. These limits are fixed points denoted by  
$$
x_+=\lim_{n\to +\infty}g^nx, \ \ \ x_-=\lim_{n\to -\infty}g^nx.
$$

\section{Results}

\paragraph{Asymptotics of weights.}

\begin{theorem}	\label{th1}
Given a point $(x,\xi)\in T^*_0M$, the weight    $\mu_{x,\xi,s}(n)$ has exponential asymptotics as $n\to \pm\infty$. More precisely, we have equivalences of weights\footnote{Recall that two positive weights $w_1(n)$ and $w_2(n)$ are equivalent  if there are uniform bounds $$c\le w_1(n)/w_2(n)\le C$$ for some positive constants $c$ and $C$.}
\begin{equation}
	\label{4}
	  \mu_{x,\xi,s}(n)\sim
	\left\{
	\begin{aligned}
	\left|\frac{\det \partial g|_{x=x_+}}{\lambda_{\min}^{2s}(x_+)}\right|^n,\ \  \xi\notin E_{+,\max}(x); \\
	\left|\frac{\det \partial g|_{x=x_+}}{\lambda_{\max}^{2s}(x_+)}\right|^n,\ \  \xi\in E_{+,\max}(x)
	\end{aligned}
	\right.\ \ \ \   \mbox{as}\ \  n\to+\infty
\end{equation}
\begin{equation}
	\label{5}
	 \mu_{x,\xi,s}(n)\sim
	\left\{
	\begin{aligned}
	\left|\frac{\det \partial g|_{x=x_-}}{\lambda_{\max}^{2s}(x_-)}\right|^n,\ \  \xi\notin E_{-,\min}(x); \\
	\left|\frac{\det \partial g|_{x=x_-}}{\lambda_{\min}^{2s}(x_-)}\right|^n,\ \  \xi\in E_{-,\min}(x)
	\end{aligned}
	\right.\ \ \ \ \mbox{as}\ \  n\to-\infty 
\end{equation}
where $\partial g|_{x=x_\pm}$ is the codifferential of $g$ at the points $x_\pm\in M^g$, while
$\lambda_{\min}(x_\pm), \lambda_{\max}(x_\pm)$ are its minimal and   maximal eigenvalues, while  $E_{+,\max}(x),E_{-,\min}(x) \subset T_{x}M$ are certain one-dimensional subspaces. 		
\end{theorem} 
The proof of Theorem~\ref{th1} relies on the following result from \cite{Hart1}, see also \cite{Bel1}: 
a germ of a hyperbolic   $C^2$-diffeomorphism $(\mathbb{R}^2,0)\to(\mathbb{R}^2,0)$ is $C^1$-conjugate to a linear diffeomorphism. 


\paragraph{Invertibility of operators with constant coefficients.}

Let
\begin{equation}
\label{6}
D=\sum_k{a_k}T^k: \ \ \ H^s(M) \longrightarrow H^s(M),
\end{equation}
be an operator with constant coefficients $a_k\in\mathbb C$. 
\begin{theorem}
	The operator  \eqref{6} is invertible if and only if there are no zeroes of the polynomial
	 $\sum_k{a_kz^k}$ in the ring  
\begin{equation}
	\label{7}
	K_s=\{z\in \mathbb{C}\;|\; r_s\le|z|\le R_s\}, 
\end{equation} 
$$
r_s=\min_{x\in M^g,j\in\{1,2\}}  (\det \partial g|_x)^{1/2}\lambda_{j}^{-s}(x),\quad
R_s=\max_{x\in M^g,j\in\{1,2\}}  (\det \partial g|_x)^{1/2}\lambda_{j}^{-s}(x)
$$
where $\lambda_{1,2}(x)$ are the two eigenvalues of $\partial g|_x$. 
\end{theorem}

\paragraph{Ellipticity of operators with variable coefficients.}
For operators with variable coefficients,  it is impossible in the general case to give explicit ellipticity conditions,
since in this case the symbol \eqref{2} is a difference operator with variable coefficients.  However, one can describe the dependence of ellipticity of a given operator on the smoothness exponent $s$ of Sobolev spaces. More precisely, the following result holds.  
\begin{theorem}
Given an operator of the form \eqref{1} associated with a diffeomorphism $g$ of Morse--Smale type, there exists an interval $I=(s_0,s_1)$, where ${-\infty\le s_0\le s_1\le +\infty}$ with the following properties:
	\begin{enumerate}
		\item[1)] if $D$ is elliptic for some $s\in I$, then it is elliptic for all $s\in I;$
		\item[2)] $D$ is not elliptic whenever $s\notin I.$
		
	\end{enumerate}
\end{theorem}


\end{document}